\documentstyle[12pt,epsfig, amsfonts,amsmath]{article}

\begin{document}

\newtheorem{thm}{Theorem}[section]
\newtheorem{thmA}[thm]{Theorem A}
\newtheorem{cor}[thm]{Corollary}
\newtheorem{prop}[thm]{Proposition}
\newtheorem{define}[thm]{Definition}
\newtheorem{rem}[thm]{Remark}
\newtheorem{example}[thm]{Example}
\newtheorem{lemma}[thm]{Lemma}
\def\theequation{\thesection.\arabic{equation}}
\newcommand{\wh}{\widehat }

\title{A new solution representation for the BBM equation in
a quarter plane and the eventual periodicity}

\author{John Meng-Kai Hong \\Department of
Mathematics\\ National Central University\\
Chung-Li 32054, Taiwan\\{\small Email:
jhong@math.ncu.edu.tw}\\
\\Jiahong Wu\\Department of Mathematics\\ Oklahoma State University
\\401 Mathematical Sciences\\ Stillwater, OK 74078\\ USA \\ {\small Email:
jiahong@math.okstate.edu}\\
\\Juan-Ming Yuan \\Department of
Applied Mathematics\\ Providence University\\
Shalu 433, Taichung Hsien, Taiwan\\{\small Email:
jmyuan@pu.edu.tw}}
\date{\today}
\maketitle

\newpage
\begin{center}
{\bf Abstract}
\end{center}
The initial- and boundary-value problem for the
Benjamin-Bona-Mahony (BBM) equation is studied in this
paper. The goal is to understand the periodic behavior
(termed as eventual periodicity) of
its solutions corresponding to periodic boundary condition
or periodic forcing. To this aim, we derive a new formula
representing solutions of this initial- and boundary-value
problem by inverting the operator $\partial_t +\alpha \partial_x
-\gamma\partial_{xxt}$ defined in the space-time quarter plane.
The eventual periodicity of the linearized BBM equation with
periodic boundary data and forcing term is established by
combining this new representation formula and the method
of stationary phase.

\vspace{.3in}
\begin{center}
{\bf AMS (MOS) Numbers}

\vspace{.1in}
35Q53, \, 35B40, \, 76B15, \, 37K40

\vspace{.3in} {\bf Keywords:}

\vspace{.1in}
BBM equation, eventual periodicity, solution representation
\end{center}

\newpage
\section{Introduction}
\setcounter{equation}{0}

This paper is concerned with the initial- and boundary-value problem
(IBVP) for the Benjamin-Bona-Mahony (BBM) equation in the quarter plane
$\{(x,t), x\ge 0, t\ge 0\}$,
\begin{equation}\label{ibvp}
\left\{
\begin{array}{ll}
\displaystyle u_t + \alpha u_x  + \beta uu_x - \gamma u_{xxt} =f,
\quad \,x\ge 0,\,t\ge0,\\
\displaystyle u(x,0) =u_0(x), \quad x\ge 0, \\
\displaystyle u(0,t) =g(t), \quad t\ge 0
\end{array}
\right.
\end{equation}
where $\gamma>0$, $\alpha$ and $\beta$ are real parameters.
Our goal is two-fold: first, to establish a new representation formula for
(\ref{ibvp}) by inverting the operator $\partial_t + \alpha\partial_x
-\gamma \partial_{xxt}$; second, to understand if the solution of (\ref{ibvp}) exhibits
certain time-periodic behavior if the boundary data $g$ and the forcing
term $f$ are periodic.

\vspace{.1in}
The study on the large-time periodic behavior is partially motivated by a
laboratory experiment involving water waves generated by a wavemaker mounted at
the end of a water channel. It is observed that if the wavemaker is oscillated
periodically, say with a long period $T_0$, it appears that in due
course, at any fixed station down the channel, the wave elevation
becomes periodic of period $T_0$. Professor Jerry L. Bona proposed the
problem of establishing this observation as a mathematically exact
fact about solutions of the suitable model equations for water
waves. One goal of this paper is to determine if the solution $u$
of (\ref{ibvp}) exhibits eventual periodicity. More precisely, we investigate
whether the difference
\begin{equation}\label{diff}
u(x,t+T_0)-u(x,t)
\end{equation}
approaches zero as $t \to \infty$ if $g$ is periodic of period $T_0$.
Since the solution of (\ref{ibvp}) grows
in time (measured in the norm of Sobolev spaces $H^k$ with $k\ge 0$)(see \cite{BW}),
the issue of eventual periodicity appears to be extremely difficult.

\vspace{.1in}
The eventual periodicity has previously been studied in several works. In \cite{BW},
Bona and Wu thoroughly investigated the large-time behavior of solutions to the BBM
equation and the KdV equation including the eventual periodicity. A formula
representing the solution of the BBM equation was derived through the Laplace
transform with respect to temporal variable $t$ and the eventual periodicity
is shown for the linearized BBM equation with zero initial data and forcing term.
It appears that
the formula derived there can not be easily extended to include a non-zero
initial data or a forcing term.  Bona, Sun and Zhang in \cite{BSZ} established
the eventual
periodicity in the context of the damped KdV equation
$$
u_t + u_x + uu_x + u_{xxx} +  u =0
$$
with small amplitude boundary data $u(0,t)=g(t)$. They were able to obtain
time-decaying bounds for solutions of this equation and the eventual periodicity
follows as a consequence. Through the Laplace transform with respect to the
spatial variable $x$, Shen, Wu and Yuan in a recent work \cite{SWY} obtained a new
solution representation formula for the KdV equation and re-established the
eventual periodicity of the linearized KdV equation. In addition, the eventual
periodicity of the full KdV equation were studied there through extensive numerical
experiments.

\vspace{.1in}
In this paper, we first derive a new solution formula for the IBVP
\begin{equation} \label{linp}
\left\{
\begin{array}{l}
\displaystyle u_t  + \alpha u_x -\gamma u_{xxt} =f, \quad x\ge 0, \,t\ge0 \\
\displaystyle u(x,0) =u_0(x), \quad x\ge 0, \\
\displaystyle u(0,t) =g(t),\quad t\ge 0.
\end{array}
\right.
\end{equation}
This explicit formula reads
\begin{eqnarray}
u(x,t) &=& g(t)e^{-\frac{x}{\sqrt{\gamma}}} + \int_0^\infty \Gamma(x-y,t) u_0(y)\,dy
\nonumber\\
&& + \int_0^t \int_0^\infty \Phi(x-y,t-\tau) \left[f(y,\tau) + \frac{\alpha}{\sqrt{\gamma}}
\,g(\tau)\, e^{-\frac{y}{\sqrt{\gamma}}} \right]\,dy d\tau,
\label{ex}
\end{eqnarray}
where $\Gamma$ and $\Phi$ are given by
\begin{eqnarray}
\Gamma(x,t) &=& \int_{-\infty}^\infty e^{ix\xi
-i\frac{\alpha\xi}{1+\gamma \xi^2} t}\, d\xi,
\nonumber\\
\Phi(x,t) &=& \int_{-\infty}^\infty \frac{1}{1+\gamma \xi^2}
e^{ix\xi-i\frac{\alpha\xi}{1+\gamma \xi^2} t}\, d\xi.
\nonumber
\end{eqnarray}
$\Gamma$ should be understood as a distribution. To obtain (\ref{ex}), we consider both
the even and odd extensions of (\ref{linp}) to the whole spatial line. By taking the
Fourier transforms of these extensions and solving the resulting equations simultaneously,
we are able to represent
\begin{equation} \label{cs}
\int_0^\infty \sin(x\xi) \, u(x,t) \,dx \quad\mbox{and}\quad
\int_0^\infty \cos(x\xi) \, u(x,t) \,dx
\end{equation}
in terms of $f$, $u_0$ and $g$. (\ref{ex}) is then established by taking the inverse Fourier
transform of the quantities in (\ref{cs}). Corollaries of (\ref{ex}) include explicit solution
formulas of (\ref{linp}) with  $\alpha=0$ and $\gamma>0$ and of (\ref{linp}) with  $\gamma=0$.
More details can be found in the second section.

\vspace{.1in}
In \cite{BCSZ} and \cite{BPS}, the IBVP (\ref{ibvp}) has been recast
as an integral equation through the inversion of the operator
$\partial_t-\gamma \partial_{xxt}$,
\begin{equation}\label{int}
u(x,t) =  u_0(x) + g(t) e^{-\frac{x}{\sqrt{\gamma}}} + \int_0^t
\int_0^\infty K(x,y) \left(\alpha u + \frac12\beta
u^2\right)(y,\tau)  \,dy\,d\tau,
\end{equation}
where the kernel function $K(x,y)$ is given by
\begin{equation}\label{kernel}
K(x,y) = \frac{1}{2\gamma} \left[e^{-\frac{x+y}{\sqrt{\gamma}}} +
\mbox{sgn}(x-y) e^{-\frac{|x-y|}{\sqrt{\gamma}}}\right].
\end{equation}
While this representation is handy in dealing with the well-posedness issue, we
find it inconvenient in studying the eventual periodicity of (\ref{ibvp}) due to
the inclusion of the linear term on the right. The new representation (\ref{ex})
allows us to show that any solution of (\ref{linp}) is eventual periodic if
$f$ and $g$ are periodic of the same period. The precise statement and its proof
are presented in the third section.

\vspace{.1in}
We mention that there is adequate theory of well-posedness and regularity on the
IBVP (\ref{ibvp}). The following theorem of Bona and Luo \cite{BL} serves
our purpose well. In the following theorem, we write ${\mathbf
R}^+=[0,\infty)$ and $C_b^k$ is
exactly like $C^k$ except that the functions
and its first $k$ derivatives are required to be bounded.
\begin{thm}
Let $I=[0,T]$ if $T$ is positive or $I=[0,\infty)$ if $T=\infty$.
Assume that $g\in C^1(I)$ and $u_0\in C_b^2({\mathbf R}^+)\cap H^2({\mathbf R}^+)$.
Then (\ref{ibvp}) is globally well-posed in the sense that there is a
unique classical solution $u\in C^1(I, C_b^\infty({\mathbf R}^+))
\cap C(I; H^2({\mathbf R}^+))$
which depends continuously on $g\in C^1(I)$
and $u_0\in C_b^2({\mathbf R}^+)\cap H^2({\mathbf R}^+)$.
\end{thm}

\vspace{.1in} A more recent work \cite{BCSZ} reduces the regularity
assumptions to $g\in C(I)$ and $u_0\in C^1_b({\mathbf R}^+)$ while
(\ref{ibvp}) still has a unique global solution in a slightly weak
sense. We shall not attempt to optimize these regularity assumptions in
this paper. The rest of this paper is divided into two sections and two
appendices.

\vspace{.3in}
\section{The inversion of the operator $\partial_t +
\alpha \partial_x -\gamma \partial_{xxt}$}
\setcounter{equation}{0} \label{inv2}

\vspace{.1in} This section explicitly solves the IBVP of the linearized
BBM equation
\begin{equation}\label{simp2}
\left\{
\begin{array}{l}
u_t  + \alpha u_x -\gamma u_{xxt} =f, \quad x\ge 0, \,t\ge0, \\
\\
u(x,0) =u_0(x), \quad x\ge 0, \\ \\
u(0,t) =g(t),\quad t\ge 0,
\end{array}
\right.
\end{equation}
where $\gamma\ge 0$ and $\alpha$ are real parameters. This amounts to
inverting the operator  $\partial_t +
\alpha \partial_x -\gamma \partial_{xxt}$ for the quarter-plane problem.
Without loss of generality, we assume $g(0)=0$.

\vspace{.1in} In the case when $\gamma>0$, we consider a new dependent variable
\begin{equation}\label{wu}
w(x,t) = u(x,t) - g(t) \,e^{-\frac{x}{\sqrt{\gamma}}},
\end{equation}
which satisfies the equations
\begin{equation}\label{weq}
\left\{
\begin{array}{l}
w_t  + \alpha w_x -\gamma w_{xxt} =\tilde{f}, \quad x\ge 0, \,t\ge 0 \\
\\
w(x,0) =u_0(x), \quad x\ge 0, \\ \\
w(0,t) =0,\quad t\ge 0.
\end{array}
\right.
\end{equation}
where
\begin{equation}\label{ftil}
\tilde{f}(x,t) =f(x,t) + \frac{\alpha}{\sqrt{\gamma}}
\,g(t)\, e^{-\frac{x}{\sqrt{\gamma}}}.
\end{equation}
The IBVP (\ref{weq}) has its boundary data equal to zero and is
solved through odd and even extensions to
the whole spatial line. The solution $u$ of (\ref{simp2}) is then
obtained by (\ref{wu}).

\vspace{.1in}
\begin{thm} \label{thm1}
Let $I=[0,T]$ if $T$ is positive or $I=[0,\infty)$ if $T=\infty$.
Let $u_0\in C^2_b({\mathbf R}^+)\cap H^2({\mathbf R}^+)$,
$g\in C^1(I)$ and $f\in C(I, L^2({\mathbf R}^+))$.
Without loss of generality, assume $u_0(0) =g(0)=0$.
Then the unique classical solution of (\ref{simp2}) can be written as
\begin{eqnarray}
u(x,t) &=& g(t)e^{-\frac{x}{\sqrt{\gamma}}} + \int_0^\infty \Gamma(x-y,t) u_0(y)\,dy
\nonumber\\
&& + \int_0^t \int_0^\infty \Phi(x-y,t-\tau) \left[f(y,\tau) + \frac{\alpha}{\sqrt{\gamma}}
\,g(\tau)\, e^{-\frac{y}{\sqrt{\gamma}}} \right]\,dy d\tau,
\label{exu}
\end{eqnarray}
where $\Gamma$ and $\Phi$ are given by
\begin{eqnarray}
\Gamma(x,t) &=& \int_{-\infty}^\infty e^{ix\xi
-i\frac{\alpha\xi}{1+\gamma \xi^2} t}\, d\xi,
\label{gam}\\
\Phi(x,t) &=& \int_{-\infty}^\infty \frac{1}{1+\gamma \xi^2}
e^{ix\xi-i\frac{\alpha\xi}{1+\gamma \xi^2} t}\, d\xi.
\label{phi}
\end{eqnarray}
$\Gamma$ in (\ref{gam}) should be understood as a distribution.
\end{thm}

\vspace{.1in}
To gain an initial understanding of the formula in this theorem, we consider
two special cases. The first is $\alpha=0$ and $f\equiv 0$.  When $\alpha=0$,
$$
\Gamma(x,t) = \int_{-\infty}^\infty e^{ix\xi} \,d\xi = \delta(x),
$$
where $\delta$ denotes the Dirac delta. Therefore, for $x\ge 0$,
$$
\int_0^\infty \Gamma(x-y,t) u_0(y)\,dy =u_0(x).
$$
\begin{cor}
The solution of (\ref{simp2}) with $\alpha=0$ and $f\equiv 0$ is given by
$$
u(x,t) =u_0(x) + g(t)e^{-\frac{x}{\sqrt{\gamma}}}.
$$
\end{cor}

\vspace{.1in}
The second special case is when $\gamma=0$ and $f\equiv 0$.
Although $\gamma=0$ in not allowed in Theorem \ref{thm1},
the solution formula for this case can still be obtained similarly as
(\ref{exu}). Instead of (\ref{wu}), one considers
$$
w(x,t)=u(x,t) -g(t)e^{-x}
$$
which solves
$$
\left\{
\begin{array}{l}
w_t  + \alpha w_x  =(\alpha g(t)-g'(t))e^{-x}, \quad x\ge 0,\, t\ge0, \\
\\
w(x,0) =u_0(x), \quad x\ge 0, \\ \\
w(0,t) =0,\quad t\ge 0.
\end{array}
\right.
$$
Then, as in Theorem \ref{thm1}, the solution for this special case is
\begin{eqnarray}
u(x,t) &=& g(t) e^{-x}  + \int_0^\infty \Gamma(x-y,t) u_0(y)\,dy
\nonumber\\
&& + \int_0^t \int_0^\infty \Phi(x-y,t-\tau)
\left[\alpha g(\tau)e^{-y}-g'(\tau)e^{-y}\right]\,dy d\tau.
\label{sp2u}
\end{eqnarray}
This representation allows us to extract a simple solution
formula for the IBVP (\ref{simp2}) with
$\gamma=0$ and $f\equiv 0$.

\begin{cor}\label{sp2}
The solution of the IBVP (\ref{simp2}) with
$\gamma=0$ and $f\equiv 0$ is given by
$$
u(x,t) =\left\{ \begin{array}{ll}
\displaystyle u_0(x-\alpha t), \quad &\mbox{if $x \ge \alpha t$}, \\
\displaystyle g\left(t-\frac{x}{\alpha}\right), \quad &\mbox{if $x < \alpha t$}.
\end{array}
\right.
$$
\end{cor}
{\it Proof of Corollary \ref{sp2}}.\quad When $\gamma=0$,
$$
\Gamma(x,t) = \int_{-\infty}^\infty \,e^{i(x-\alpha t)\xi}\,d\xi
= \delta(x-\alpha t),\qquad \Phi(x,t) = \delta(x-\alpha t).
$$
If $x-\alpha t\ge 0$, then
\begin{equation}\label{tm1}
\int_0^\infty \Gamma(x-y,t) u_0(y)\,dy =\int_0^\infty \delta(x-\alpha t -y)
u_0(y) \,dy =u_0(x-\alpha t)
\end{equation}
and
\begin{eqnarray}
&& \int_0^t \int_0^\infty \Phi(x-y,t-\tau)
\left[\alpha g(\tau)e^{-y}-g'(\tau)e^{-y}\right]\,dy d\tau \nonumber \\
&& \qquad = \int_0^t \left[\alpha g(\tau)-g'(\tau)\right] \int_0^\infty
\delta(x-\alpha t+\alpha \tau -y)\, e^{-y} \,dy \,d\tau \nonumber \\
&& \qquad =\int_0^t \left[\alpha g(\tau)-g'(\tau)\right]
e^{-(x-\alpha t+\alpha \tau)} \,d\tau \nonumber \\
&& \qquad=e^{-x + \alpha t} \left[\int_0^t \alpha g(\tau)\,e^{-\alpha \tau}\,d\tau
-\int_0^t g'(\tau)\,e^{-\alpha \tau}\,d\tau \right]  \nonumber \\
&& \qquad=-g(t) e^{-x}. \label{tm2}
\end{eqnarray}
Inserting (\ref{tm1}) and (\ref{tm2}) in (\ref{sp2u}) yields
$$
u(x,t) =u_0(x-\alpha t) \qquad \mbox{if $x-\alpha t\ge 0$}.
$$
If $x-\alpha t <0$, then
$$
\int_0^\infty \delta(x-\alpha t -y)
u_0(y) \,dy = 0.
$$
and
$$
\int_0^\infty
\delta(x-\alpha t+\alpha \tau -y)\, e^{-y} \,dy = \left\{
\begin{array}{ll}
\displaystyle e^{-(x-\alpha t+\alpha \tau)}, \quad &\mbox{if $\tau \ge t-\alpha^{-1} x$}, \\
\displaystyle 0, \quad &\mbox{otherwise}.
\end{array}
\right.
$$
Thus,
\begin{eqnarray}
&& \int_0^t \int_0^\infty \Phi(x-y,t-\tau)
\left[\alpha g(\tau)e^{-y}-g'(\tau)e^{-y}\right]\,dy d\tau \nonumber \\
&& \qquad =\int_{t-\alpha^{-1} x}^t e^{-(x-\alpha t+\alpha \tau)}
\left[\alpha g(\tau)-g'(\tau)\right]\,d\tau
\nonumber\\
&& \qquad = -g(t) e^{-x} + g(t-\alpha^{-1}\,x). \nonumber
\end{eqnarray}
Therefore, for $x<\alpha t$,
$$
u(x,t) = g(t-\alpha^{-1}\,x).
$$
This completes the proof of Corollary \ref{sp2}.

\vspace{.15in}
\noindent {\it Proof of Theorem \ref{thm1}}.\quad The major idea is to
extend the equation of (\ref{weq}) from the half line $\{x: x>0\}$
to the whole line $x\in {\mathbf R}$ so that the method of Fourier
transforms can be employed. Both odd and even extensions will be considered.
The rest of the proof is divided into four major steps.

\vspace{.1in} {\it Step 1. Odd extension}.
Denote by $W$ the odd extension of $w$ in $x$, namely
$$
W(x,t) =
\left\{
\begin{array}{ll} w(x,t) \quad &\mbox{if $x\ge 0$,}\\
-w(-x,t) \quad &\mbox{if $x<0$.}
\end{array}
\right.
$$
If $W_0$ and $F$ are the odd extensions of $w_0$ and $\tilde{f}$, respectively,
then $W$ solves the following initial-value problem
\begin{equation}\label{ww}
\left\{
\begin{array}{l}
W_t  + \alpha\, \mbox{sgn}(x)W_x -\gamma W_{xxt} =F, \quad x \in {\mathbf R}, \,t>0 \\
\\
W(x,0) =W_0(x), \quad x \in {\mathbf R}.
\end{array}
\right.
\end{equation}
Let $\widehat{W}$ denote the Fourier transform of $W$, namely
$$
\widehat{W}(\xi,t)= {\cal F}(W)(\xi,t) = \int_{-\infty}^\infty
e^{- i\xi x} W(x,t) dx.
$$
 After applying a basic property of the Fourier transform, we obtain
\begin{equation}\label{fou}
(1+ \gamma \xi^2) \,\partial_t \widehat{W}(\xi,t) + \alpha \,{\cal F}\left(\mbox{sgn}(x)W_x\right)
\,=\,\widehat{F}(\xi,t).
\end{equation}
According to the definition of $W$, we have
\begin{equation}\label{wh}
\widehat{W}(\xi,t) =\int_0^\infty (e^{-ix\xi}-e^{ix\xi}) w(x,t) dx
=-2i \int_0^\infty \sin(x\xi) w(x,t) dx.
\end{equation}
In addition,
\begin{eqnarray}
{\cal F}\left(\mbox{sgn}(x)W_x \right)(\xi,t) &=& \int_{-\infty}^0 e^{-i\xi
x} \mbox{sgn}(x)W_x \,dx + \int_0^\infty e^{-i\xi
x} \mbox{sgn}(x)W_x \,dx\nonumber\\
&=& \int_{-\infty}^0 e^{-i\xi x} (-w_x(-x,t)) dx + \int_0^\infty e^{-i\xi x}
w_x(x,t) dx. \nonumber
\end{eqnarray}
Making the substitution $y=-x$ and integrating by parts yield
\begin{eqnarray}
{\cal F}\left(\mbox{sgn}(x)W_x \right)(\xi,t)
&=& -\int_0^\infty e^{iy\xi} w_y(y,t)\,dy + \int_0^\infty e^{-i\xi x}
w_x(x,t) dx. \nonumber \\
&=&i\xi \int_0^\infty (e^{ix\xi} + e^{-ix\xi})\, w(x,t) \,dx
\nonumber\\
&=& 2i \xi \int_0^\infty \cos(x\xi)\, w(x,t) \,dx.\label{alt}
\end{eqnarray}
We can also write $\widehat{F}$ in terms of $\tilde{f}$ as
\begin{equation}\label{ft}
\widehat{F}(\xi,t) = -2i \int_0^\infty \sin(x\xi)\, \tilde{f}(x,t) \,dx
\end{equation}
Inserting (\ref{wh}),(\ref{alt}) and (\ref{ft}) in (\ref{fou}), we obtain
\begin{equation}\label{eq1}
\partial_t X(\xi,t) -\beta(\xi) Y(\xi,t) \,=\,
h_1(\xi,t),
\end{equation}
where
$$
X(\xi,t) = \int_0^\infty  \sin(x\xi)\, w(x,t) \,dx, \quad
Y(\xi,t) = \int_0^\infty  \cos(x\xi)\, w(x,t) \,dx,
$$
$$
\beta(\xi) = \frac{\alpha \xi}{1+\gamma \xi^2} \quad \mbox{and}\quad h_1(\xi,t)
=\frac{1}{1+\gamma \xi^2}\,\int_0^\infty \sin(x\xi)\, \tilde{f}(x,t) \,dx.
$$

\vspace{.1in}
{\it Step 2. Even extension}. Denote by $V(x,t)$ the even extension of $w$, namely
$$
V(x,t) =
\left\{
\begin{array}{ll} w(x,t) \quad &\mbox{if $x\ge 0$,}\\
w(-x,t) \quad &\mbox{if $x<0$.}
\end{array}
\right.
$$
It can be easily verified that $V$ satisfies
\begin{equation}\label{vw}
\left\{
\begin{array}{l}
V_t  + \alpha\, \mbox{sgn}(x)V_x -\gamma V_{xxt} =H, \quad x \in {\mathbf R}, \,t>0 \\
\\
V(x,0) =V_0(x), \quad x \in {\mathbf R} ,
\end{array}
\right.
\end{equation}
where $H$ and $V_0$ are the even extensions of $\tilde{f}$ and $w_0$, respectively.
As in {\it Step 1}, we have
\begin{equation}\label{v1}
\widehat{V}(\xi,t) = 2 \int_0^\infty \cos(x\xi)\,w(x,t) \,dx,
\end{equation}
and
\begin{equation}\label{v2}
\widehat{H}(\xi,t) = 2 \int_0^\infty \cos(x\xi)\, \tilde{f}(x,t) \,dx
\end{equation}
Furthermore,
\begin{eqnarray}
{\cal F}\left(\mbox{sgn}(x)V_x \right)(\xi,t) &=& \int_{-\infty}^0 e^{-i\xi
x} \mbox{sgn}(x)V_x \,dx + \int_0^\infty e^{-i\xi
x} \mbox{sgn}(x)V_x \,dx\nonumber\\
&=& \int_{-\infty}^0 e^{-i\xi x} w_x(-x,t) dx + \int_0^\infty e^{-i\xi x}
w_x(x,t) dx. \nonumber
\end{eqnarray}
Making the substitution $y=-x$ in the first integral and integrating by parts leads to
\begin{eqnarray}
{\cal F}\left(\mbox{sgn}(x)V_x \right)(\xi,t) &=&
 \int_0^\infty e^{i\xi y}  w_y(y,t)\,dy + \int_0^\infty e^{-i\xi x}
 w_x(x,t)\,dx
\nonumber \\
&=& -i\xi \int_0^\infty e^{i\xi x}  w(x,t) \,dx
+ i\xi \int_0^\infty e^{-i\xi x} w(x,t) \,dx \nonumber \\
&=& 2\xi \int_0^\infty \sin(x\xi)\,w(x,t)\,dx. \label{v3}
\end{eqnarray}
Taking the Fourier transform of (\ref{vw}) and applying (\ref{v1}),(\ref{v2})
and (\ref{v3}), we obtain
\begin{equation}\label{eq2}
\partial_t Y(\xi,t) + \beta X(\xi,t) = h_2(\xi,t),
\end{equation}
where
$$
h_2(\xi,t)= \frac{1}{1+\gamma\xi^2} \,
\int_0^\infty \cos(x\xi) \, \tilde{f}(x,t)\,dx.
$$

\vspace{.1in}
{\it Step 3. Solving for $X(\xi,t)$ and $Y(\xi,t)$}. Solving the
system of (\ref{eq1}) and (\ref{eq2}), we find
$$
X(\xi,t) = \int_0^\infty \sin(x\xi + \beta t) \,u_0(x)dx
+ \frac{1}{1+\gamma \xi^2}\int_0^t \int_0^\infty \sin(x\xi + \beta(t-\tau)) \,
\tilde{f}(x,\tau) dx d\tau,
$$
$$
Y(\xi,t) = \int_0^\infty \cos(x\xi + \beta t) \,u_0(x)dx
+ \frac{1}{1+\gamma \xi^2}\int_0^t \int_0^\infty \cos(x\xi + \beta(t-\tau)) \,
\tilde{f}(x,\tau) dx d\tau.
$$
where $\tilde{f}$ is defined in (\ref{ftil}). We leave the details of derivation
in Appendix A.

\vspace{.1in}
{\it Step 4. Finding $w(x,t)$ through the inverse Fourier transform}.
To find $w(x,t)$, we first notice that
$$
\int_0^\infty e^{-i y \xi} w(y,t) dy = Y(\xi,t) -i X(\xi,t).
$$
Applying the formulas in the previous step, we have
\begin{eqnarray}
\int_0^\infty e^{-i y \xi} w(y,t) dy &=& e^{-i\beta t}\,\int_0^\infty e^{-iy\xi} u_0(y) dy
\nonumber\\
&& +  \frac{1}{1+\gamma \xi^2}\int_0^t e^{-i \beta (t-\tau)}
\int_0^\infty e^{-i y\xi}\,
\tilde{f}(y,\tau) dy d\tau. \label{wfo}
\end{eqnarray}
We shall now establish a theorem that allows us to obtain $w(x,t)$ by
taking the inverse Fourier transform of (\ref{wfo}).
\begin{thm} \label{invf}
Fix $t>0$. If $u(x,t) \in L^2({\mathbf R}^+)$, then, for any $x\ge 0$,
\begin{equation}\label{ff}
\int_{-\infty}^\infty e^{ix\xi} \int_0^\infty e^{-iy\xi} u(y,t)
\,dy\,d\xi =u(x,t).
\end{equation}
\end{thm}
{\it Proof of Theorem \ref{invf}}. Recall that if $ g_\epsilon(\xi)
= \exp{(-\epsilon \pi \xi^2)}$, then
$$
\widehat{g}_\epsilon (x)= \frac{1}{\epsilon^{1/2}}
\exp{\left(-\frac{\pi x^2}{\epsilon}\right)}.
$$
In addition, for any $f\in L^p({\mathbf R})$ with $1\le p<\infty$,
$$
\widehat{g}_\epsilon \ast f \,\,\to \,\, f \quad \mbox{in} \quad
L^p({\mathbf R}).
$$
These basic facts can be found in Lieb and Loss
\cite{LL}. Now, consider
$$
\int_{-\infty}^\infty e^{ix\xi} g_\epsilon(\xi) \int_0^\infty
e^{-iy\xi} u(y,t) \,dy\,d\xi.
$$
According to Lemma \ref{l2l} below, if $u(y,t) \in L^2({\mathbf
R}^+)$, then
$$
P(u)(\xi,t) \equiv \int_0^\infty e^{-iy\xi} u(y,t) \,dy \in
L^2({\mathbf R}).
$$
Since $g_\epsilon(\xi) \to 1$ as $\epsilon\to 0$, the dominated
convergence theorem implies
$$
g_\epsilon P(u) \,\,\to \,\, P(u) \quad\mbox{in}\quad
L^2({\mathbf{R}}).
$$
Then,
\begin{eqnarray}
\int_{-\infty}^\infty e^{ix\xi}\, g_\epsilon(\xi) \int_0^\infty
e^{-iy\xi} u(y,t) \,dy\,d\xi &=& \int_0^\infty u(y,t)
\int_{-\infty}^\infty e^{-i(y-x)\xi} g_\epsilon(\xi)\,d\xi\,dy
\nonumber \\
&=&\int_0^\infty u(y,t) \widehat{g}_\epsilon(y-x) dy. \label{ff1}
\end{eqnarray}
As in the proof of Theorem 2.16 of \cite{LL}, we can prove
$$
\int_0^\infty u(y,t) \widehat{g}_\epsilon(y-x) dy \,\,\to\,\, u(x,t)
\quad\ \mbox{in}\quad L^2({\mathbf{R}}^+).
$$
Letting $\epsilon\to 0$ in (\ref{ff1}) yields (\ref{ff}).
This proves Theorem \ref{invf}.

\begin{lemma} \label{l2l}
If $u(y,t) \in L^2({\mathbf R}^+)$, then
$$
P(u)(\xi,t) \equiv \int_0^\infty e^{-iy\xi} u(y,t) \,dy \in
L^2({\mathbf R}).
$$
\end{lemma}
{\it Proof of Lemma \ref{l2l}}. For any $\epsilon>0$,
\begin{eqnarray}
\int_{-\infty}^\infty |P(u)|^2(\xi,t) g_\epsilon(\xi) \,d\xi &=&
\int_{-\infty}^\infty P(u)(\xi,t)\overline{P(u)}(\xi,t)
g_\epsilon(\xi) \,d\xi \nonumber \\
&=& \int_{-\infty}^\infty g_\epsilon(\xi) \,\int_0^\infty
e^{-ix\xi}u(x,t) dx \,\int_0^\infty e^{iy\xi}u(y,t) dy \,d\xi \nonumber \\
&=& \int_0^\infty \int_0^\infty u(x,t) u(y,t) \int_{-\infty}^\infty
e^{-i(x-y)\xi} g_\epsilon(\xi) \,d\xi dx dy \nonumber \\
&=& \int_0^\infty u(x,t) \int_0^\infty \widehat{g}_\epsilon(x-y)
u(y,t) dy \,dx \nonumber
\end{eqnarray}
As $\epsilon\to 0$,
$$
\int_0^\infty \widehat{g}_\epsilon(x-y) u(y,t) dy \,\, \to
\,\,u(x,t) \quad\mbox{in}\quad L^2({\mathbf{R}}^+)
$$
and $u(x,t) \in L^2({\mathbf R}^+)$ implies $\int_{-\infty}^\infty
|P(u)|^2(\xi,t) g_\epsilon(\xi) \,d\xi$ is bounded uniformly. Since
$$
g_\epsilon(\xi) = \exp(-\epsilon \pi |\xi|^2) \quad\mbox{$\to 1$ as
$\epsilon \to 0$},
$$
we obtain by applying the monotone convergence theorem that
$$
\int_{-\infty}^\infty |P(u)|^2(\xi,t) \,d\xi = \int_0^\infty
|u(x,t)|^2 \,dx.
$$
This proves Lemma \ref{l2l}.

\vspace{.1in}
Taking the inverse Fourier transform (denoted by ${\cal F}^{-1}$)
of (\ref{wfo}) and applying Theorem \ref{invf} and the basic property
${\cal F}^{-1}(fg)={\cal F}^{-1}(f)\ast {\cal F}^{-1}(g)$, we obtain
$$
w(x,t) = \int_0^\infty \Gamma(x-y,t) u_0(y)\,dy
+ \int_0^t \int_0^\infty \Phi(x-y,t-\tau) \tilde{f}(y,\tau) \,dy d\tau,
$$
where, noticing $\beta =\frac{\alpha\xi}{1+\gamma \xi^2}$,
\begin{eqnarray}
\Gamma(x,t) &=& \int_{-\infty}^\infty e^{ix\xi
-i\frac{\alpha\xi}{1+\gamma \xi^2} t}\, d\xi,
\nonumber\\
\Phi(x,t) &=& \int_{-\infty}^\infty \frac{1}{1+\gamma \xi^2}
e^{ix\xi-i\frac{\alpha\xi}{1+\gamma \xi^2} t}\, d\xi.
\nonumber
\end{eqnarray}
Therefore, by (\ref{weq}),
\begin{eqnarray}
u(x,t) &=& w(x,t) + g(t)e^{-\frac{x}{\sqrt{\gamma}}}
\nonumber\\
&=& g(t)e^{-\frac{x}{\sqrt{\gamma}}} + \int_0^\infty \Gamma(x-y,t) u_0(y)\,dy
\nonumber\\
&& + \int_0^t \int_0^\infty \Phi(x-y,t-\tau) \left[f(y,\tau)
+ \frac{\alpha}{\sqrt{\gamma}}
\,g(\tau)\, e^{-\frac{y}{\sqrt{\gamma}}} \right]\,dy d\tau.
\nonumber
\end{eqnarray}
This completes the proof of Theorem \ref{thm1}.

\vspace{.3in}
\section{Eventual periodicity}
\setcounter{equation}{0}
\label{evenp}

This section studies the eventual periodicity of solutions to the IBVP
for the linearized BBM equation
\begin{equation}\label{libvp}
\left\{
\begin{array}{ll}
\displaystyle u_t + \alpha u_x  - \gamma u_{xxt} =f(x,t),\quad \,x\ge 0,\,t\ge0,\\
\displaystyle u(x,0) =u_0(x), \quad x\ge 0, \\
\displaystyle u(0,t) =g(t), \quad t\ge 0
\end{array}
\right.
\end{equation}
where $g$ and $f$ are assumed to be periodic in $t$.   We prove the following
theorem.

\begin{thm} \label{lincase}
Let $\gamma>0$ and $\alpha$ be real parameters.
Let $I=[0,T]$ if $T$ is positive or $I=[0,\infty)$ if $T=\infty$.
Let $u_0\in C_b^2({\mathbf R}^+)$ and $u_0' \in L^1({\mathbf R}^+)$. Let $g \in C(I)$
with $g(0)=0$ and $f\in C(I, L^1({\mathbf R}^+))$. Assume $g$ and $f$ are
periodic of period $T_0$ in $t$, namely, for all $t\ge 0$,
\begin{equation}\label{pc}
g(t+T_0)=g(t) \quad\mbox{and}\quad f(x,t+T_0)=f(t).
\end{equation}
Then, for any fixed $x>0$, the solution $u$
of (\ref{libvp}) satisfies
\begin{equation}\label{even}
\lim_{t\to \infty} \left(u(x,t+T_0)-u(x,t) \right) =0.
\end{equation}
That is, $u$ is eventually periodic of period $T_0$.
\end{thm}

\begin{rem}
When $\gamma=0$, the eventual periodicity is easily obtained from
the explicit formula in Corollary \ref{sp2}.
\end{rem}

\vspace{.1in}
\noindent {\it Proof of Theorem \ref{lincase}}. \quad
Consider the new function
$$
v(x,t) =u(x,t)-u_0(x),
$$
which satisfies
\begin{equation}\label{veq}
\left\{
\begin{array}{ll}
\displaystyle v_t + \alpha v_x  - \gamma v_{xxt}
=f(x,t)-\alpha \, u_0'(x),\quad \,x\ge 0,\,t\ge0,\\
\displaystyle v(x,0) =0, \quad x\ge 0, \\
\displaystyle v(0,t) =g(t), \quad t\ge 0
\end{array}
\right.
\end{equation}
Applying the explicit formula to (\ref{veq}) gives
\begin{eqnarray}
&& v(x,t) = g(t)e^{-\frac{x}{\sqrt{\gamma}}} \nonumber \\
&& \qquad +\int_0^t \int_0^\infty \Phi(x-y,t-\tau) \left[f(y,\tau)-\alpha u_0'(y)
+ \frac{\alpha}{\sqrt{\gamma}}
\,g(\tau)\, e^{-\frac{y}{\sqrt{\gamma}}} \right]\,dy d\tau,\nonumber
\end{eqnarray}
where $\Phi$ is given by (\ref{phi}).
Noticing the conditions in (\ref{pc}),  we obtain after making a substitution,
\begin{eqnarray}
&& u(x,t+T_0) - u(x,t) = v(x,t+T_0) - v(x,t) \nonumber \\
&& \qquad  =\int_0^{t+T_0} \int_0^\infty \Phi(x-y,t+T_0-\tau) \left[f(y,\tau)-\alpha u_0'(y)+ \frac{\alpha}{\sqrt{\gamma}}
\,g(\tau)\, e^{-\frac{y}{\sqrt{\gamma}}} \right]\,dy d\tau \nonumber \\
&& \qquad \quad -\int_0^t \int_0^\infty \Phi(x-y,t-\tau) \left[f(y,\tau)-\alpha u_0'(y)+ \frac{\alpha}{\sqrt{\gamma}}
\,g(\tau)\, e^{-\frac{y}{\sqrt{\gamma}}} \right]\,dy d\tau \nonumber \\
&& \qquad =\int_{-T_0}^0 \int_0^\infty \Phi(x-y,t-\tau) \left[f(y,\tau)-\alpha u_0'(y)+ \frac{\alpha}{\sqrt{\gamma}}
\,g(\tau)\, e^{-\frac{y}{\sqrt{\gamma}}} \right]\,dy d\tau.\label{ud}
\end{eqnarray}
To show (\ref{even}), we first show that, for any $\epsilon>0$, there is $K>0$ such
that
$$
\left| \Phi(x,t) \right| = \left| \int_{-\infty}^\infty \frac{1}{1+\gamma \xi^2}\,e^{ix\xi
-i\frac{\alpha\xi}{1+\gamma \xi^2} t}\, d\xi \right| <\epsilon \quad\mbox{when $t>K$}.
$$
First, we choose $M =M(\epsilon)>0$ such that
$$
\int_{-\infty}^{-M} \frac{1}{1+\gamma \xi^2}
\, d\xi\, + \int_{M}^\infty \frac{1}{1+\gamma \xi^2}
\, d\xi \, <\, \frac{\epsilon}{2}.
$$
Next, we apply the method of stationary phase to show the following asymptotics.
\begin{prop}\label{asm}
For any fixed $M>0$,
$$
\int_{-M}^M \frac{1}{1+\gamma \xi^2}\,e^{ix\xi
-i\frac{\alpha\xi}{1+\gamma \xi^2} t}\, d\xi  =O \left(\frac{1}{\sqrt{t}}\right)
$$
as $t\to \infty$. This large-time asymptotics is uniform in $x\in {\mathbf R}^+$.
\end{prop}

We remark that the result of
this proposition does not necessarily hold for
$$
\int_{-\infty}^\infty\frac{1}{1+\gamma \xi^2}\,e^{ix\xi
-i\frac{\alpha\xi}{1+\gamma \xi^2} t}\, d\xi
$$
because one of the conditions in the method of stationary phase is violated. More
details on this point will be provided in Appendix B.
The following estimate is a special consequence of this proposition.
\begin{cor}\label{mms}
There exists $K=K(M)$ such that
$$
\left| \int_{-M}^M \frac{1}{1+\gamma \xi^2}\,e^{ix\xi
-i\frac{\alpha\xi}{1+\gamma \xi^2} t}\, d\xi \right| \, <\, \frac{\epsilon}{2}
$$
whenever $t>K$.
\end{cor}

We resume the proof of Theorem \ref{lincase}. It then follows
from (\ref{ud}) that, for $t>K$,
\begin{eqnarray}
&& \left|u(x,t+T) - u(x,t)\right| \nonumber \\
&& \qquad \le \epsilon \,\left[\int_{-T_0}^0
\int_0^\infty |f(y,\tau)|dy d\tau + \alpha T \int_0^\infty |u'_0(y)|dy
+ \alpha\,\int_{-T_0}^0 g(\tau)\,d\tau \right].\nonumber
\end{eqnarray}
This completes the proof of Theorem \ref{lincase}.

\vspace{.1in}
\noindent{\it Proof of Proposition \ref{asm}}. \quad We apply the method of stationary
phase. To do so, we split the integral into four parts
\begin{eqnarray}
\int_{-M}^M \frac{1}{1+\gamma \xi^2}\,e^{ix\xi
-i\frac{\alpha\xi}{1+\gamma \xi^2} t}\, d\xi &=& I_1 +I_2 +I_3+I_4 \nonumber \\
& \equiv & \int_{-M}^{-\frac{1}{\sqrt{\gamma}}}
\,+\, \int_{-\frac{1}{\sqrt{\gamma}}}^0 \,+\, \int_0^{\frac{1}{\sqrt{\gamma}}} \,+\,
\int_{\frac{1}{\sqrt{\gamma}}}^M \,.\nonumber
\end{eqnarray}
Note that $\pm 1/\sqrt{\gamma}$ are the zero points of the derivative of
$p(\xi)=\alpha\xi/(1+\gamma \xi^2)$ and $p'(\xi)$ is nonzero for $\xi$ in
each of these intervals. It suffices to consider $I_2$ and $I_4$. Without loss
of generality, we assume $\alpha>0$.
Direct applications of the method of
stationary phase concludes that
\begin{equation}\label{i2}
I_2  \quad\sim\quad
\frac{\sqrt{\pi}\,e^{-\frac{\pi i }{4}}\,
e^{\frac{i}{\sqrt{\gamma}}(-x+\frac{\alpha t}{2})}}{2\,\sqrt{2\alpha}
\sqrt[4]{\gamma}\,
\sqrt{t}},
\end{equation}
\begin{equation}\label{i4}
I_4 \quad \sim \quad \frac{\sqrt{\pi}\,e^{\frac{\pi i }{4}}\,
e^{\frac{i}{\sqrt{\gamma}}(x-\frac{\alpha t}{2})}}{2\,\sqrt{2\alpha}
\sqrt[4]{\gamma}\,
\sqrt{t}}.
\end{equation}
For readers' convenience, this method is recalled in Appendix B and the details
leading to these estimates can also be found there.

\vspace{.3in}
\noindent {\bf Acknowledgements}

\vspace{.1in}
The authors thank Professor Jerry Bona for initiating this problem and for his
guidance. The assistance of three students, Hui-Sheng Chen, Chieh-Peng Chuang
and Man-Meng Io, is also greatly appreciated.

\vspace{.1in}
Hong is partially supported by a NSC grant. Yuan is partially
supported by the NSC grant \#96-2115-M126-001 and
he thanks Professors Goong Chen, I-Liang Chern and Jyh-Hao Lee
for their encouragement. Both Hong and Yuan thank the Department of
Mathematics at Oklahoma State University for support and hospitality.
Wu thanks the Department of Mathematics at
National Central University (Taiwan) and the Department of Applied
Mathematics at National Chiao-Tung University (Taiwan) for support
and hospitality.

\vspace{.4in} \noindent{\large \bf Appendix A}

\renewcommand{\thesection}{\Alph{section}}
\renewcommand{\theequation}{\Alph{section}.\arabic{equation}}
\setcounter{section}{1} \setcounter{equation}{0} \setcounter{thm}{0}

\vspace{.15in} \noindent
This appendix provides the details of solving the system of ODEs
(\ref{eq1}) and (\ref{eq2}). Consider the general nonhomogeneous linear systems
$$
\frac{d}{dt} {\bf x}(t) = P(t) {\bf x}(t) + {\bf g}(t), \qquad {\bf x}(t_0) ={\bf x}_0
$$
where ${\bf x}\in {\mathbf R}^n$, $P\in {\mathbf R}^{n\times n}$
and ${\bf g}\in {\mathbf R}^n$. By variation of parameters, its solution can
be written as
\begin{equation}\label{form}
{\bf x}(t) = \Psi(t) \Psi^{-1}(t_0) {\bf x}_0 + \Psi(t) \int_{t_0}^t
\Psi^{-1}(s) {\bf g}(s) \,ds
\end{equation}
where $\Psi$ denotes a fundamental matrix of the homogeneous system
$$
\frac{d}{dt} {\bf x}(t) = P(t) {\bf x}(t).
$$
Since the system of equations we are solving can be written as
$$
\partial_t \left[ \begin{array}{c} X(\xi,t)\\ Y(\xi,t) \end{array} \right]
= \left[\begin{array}{cc} 0& \beta \\ -\beta &0 \end{array} \right]
\left[ \begin{array}{c} X(\xi,t)\\ Y(\xi,t) \end{array} \right]
+ \left[ \begin{array}{c} h_1(\xi,t)\\ h_2(\xi,t) \end{array} \right]
$$
and a fundamental matrix of the corresponding homogenous system
is given by
$$
\Psi(t) = \left[ \begin{array}{cc} e^{i\beta t} & e^{-i\beta t}\\
i\,e^{i\beta t} & -i\,e^{-i\beta t} \end{array} \right],
$$
we apply (\ref{form}) to obtain that
\begin{equation}\label{conc}
\left[ \begin{array}{c} X(\xi,t)\\ Y(\xi,t) \end{array} \right]
= \Psi(t) \Psi^{-1}(0) \left[ \begin{array}{c} X(\xi,0)\\ Y(\xi,0)
\end{array} \right] + \Psi(t)\int_0^t \Psi^{-1}(\tau)
\left[ \begin{array}{c} h_1(\xi,\tau)\\ h_2(\xi,\tau) \end{array}
\right] \, d\tau.
\end{equation}
Inserting
$$
\Psi^{-1}(t) = \frac12\left[ \begin{array}{cc} e^{-i\beta t} & -i e^{-i\beta t}\\
e^{i\beta t} & i e^{i\beta t} \end{array} \right], \quad
\left[ \begin{array}{c} X(\xi,0)\\ Y(\xi,0)
\end{array} \right]  = \left[ \begin{array}{c}  \displaystyle
\int_0^\infty \sin(x\xi) u_0(x)\,dx \\ \displaystyle\int_0^\infty \cos(x\xi)
u_0(x)\,dx \end{array}\right]
$$
and
$$
\left[ \begin{array}{c} h_1(\xi,t)\\ h_2(\xi,t) \end{array}
\right] =  \left[ \begin{array}{c}  \displaystyle
\frac{1}{1+\gamma \xi^2}\int_0^\infty \sin(x\xi) \tilde{f}(x,t)\,dx
\\ \displaystyle \frac{1}{1+\gamma \xi^2} \int_0^\infty \cos(x\xi)
\tilde{f}(x,t)\,dx \end{array}\right]
$$
in (\ref{conc}), we find after some simplification
$$
X(\xi,t) = \int_0^\infty \sin(x\xi + \beta t) \,u_0(x)dx
+ \frac{1}{1+\gamma \xi^2}\int_0^t \int_0^\infty \sin(x\xi + \beta(t-\tau)) \,
\tilde{f}(x,\tau) dx d\tau,
$$
$$
Y(\xi,t) = \int_0^\infty \cos(x\xi + \beta t) \,u_0(x)dx
+ \frac{1}{1+\gamma \xi^2}\int_0^t \int_0^\infty \cos(x\xi + \beta(t-\tau)) \,
\tilde{f}(x,\tau) dx d\tau.
$$

\vspace{.4in} \noindent{\large \bf Appendix B}

\renewcommand{\thesection}{\Alph{section}}
\renewcommand{\theequation}{\Alph{section}.\arabic{equation}}
\setcounter{section}{2} \setcounter{equation}{0} \setcounter{thm}{0}

\vspace{.15in} \noindent
This appendix offers an expanded commentary on the asymptotic analysis
of the oscillatory integrals arising in section \ref{evenp}. These
analysis relies upon standard results in
the theory of stationary phase, e.q. Theorem 13.1 in F. Olver's book
\cite{Ol}. For readers' convenience, we first recall
this theory here.

\vspace{.1in} Suppose that in the integral
$$
I(t) = \int_a^b e^{itp(y)} \, q(y)\, dy
$$
the limits $a$ and $b$ are independent of $t$, $a$ being finite and
$b (>a)$ finite or infinite. The functions $p(y)$ and $q(y)$ are
independent of $t$, $p(y)$ being real and $q(y)$ either real or
complex. We also assume that the only point at which $p'(y)$
vanishes is $a$. Without loss of generality, both $t$ and $p'(y)$
are taken to be positive; cases in which one of them is negative can
be handled by changing the sign of $i$ throughout. We require
\begin{enumerate}
\item[(i)] In $(a,b)$, the functions $p'(y)$ and $q(y)$ are
continuous, $p'(y)>0$, and $p''(y)$ and $q'(y)$ have at most a
finite number of discontinuities and infinities.
\item[(ii)] As $y\to a+$,
\begin{equation}\label{ytoa}
p(y)-p(a) \,\sim\, P(y-a)^\mu,\qquad q(y)\,\sim\,
Q(y-a)^{\lambda-1},
\end{equation}
the first of these relations being differentiable. Here $P$, $\mu$
and $\lambda$ are positive constants, and $Q$ is a real or complex
constant.
\item[(iii)] For each $\epsilon\in (0, b-a)$,
$$
{\cal V}_{a+\epsilon, b}\Big\{\frac{q(y)}{p'(y)}\Big\}\equiv
\int_{a+\epsilon}^b \left|\Big(\frac{q(y)}{p'(y)}\Big)'\right|\, dy
<\infty.
$$
\item[(iv)] As $t\to b-$, the limit of $q(y)/p'(y)$ is finite, and
this limit is zero if $p(b)=\infty$.
\end{enumerate}
With these conditions, the nature of asymptotic approximation to
$I(t)$ for large $t$ depends on the sign of $\lambda-\mu$. In the
case $\lambda<\mu$, we have the following theorem.

\begin{thm}\label{phase}
In addition to the above conditions, assume that $\lambda<\mu$, the
first of (\ref{ytoa}) is twice differentiable, and the second of
(\ref{ytoa}) is differentiable, then
$$
I(t) \quad \sim \quad e^{\lambda \pi i/(2\mu)}\, \frac{Q}{\mu}
\Gamma\Big(\frac{\lambda}{\mu}\Big) \,
\frac{e^{itp(a)}}{(Pt)^{\lambda/\mu}}
$$
as $t\to\infty$.
\end{thm}

\vspace{.15in} We now provide the details leading to (\ref{i2}).
It suffices to check the conditions of Theorem \ref{phase}. Setting
$$
a= -\frac{1}{\sqrt{\gamma}}, \quad b=0,\quad p(\xi) = \frac{\alpha \xi}{1+\gamma \xi^2}
\quad\mbox{and}\quad q(\xi) =\frac{e^{i\xi x}}{1+\gamma \xi^2},
$$
we have
\begin{enumerate}
\item[(i)] $p$, $p'$, $p''$, $q$ and $q'$ are all continuous
in $(-1/\sqrt{\gamma},0)$, and $p'(\xi)>0$.
\item[(ii)] As $\xi \to -\frac{1}{\sqrt{\gamma}}+$,
$$
p(\xi)-p\Big(-\frac{1}{\sqrt{\gamma}}\Big)  \,\,\sim \,\,
\frac{\alpha\sqrt{\gamma}}{2}
\Big(\xi+ \frac{1}{\sqrt{\gamma}}\Big)^2,\qquad
q\Big(-\frac{1}{\sqrt{\gamma}}\Big)\,\,\sim \,\,\frac12
e^{-i\frac{1}{\sqrt{\gamma}}x}.
$$
That is, $P=\frac{\alpha\sqrt{\gamma}}{2}$, $\mu=2$, $Q=\frac12
e^{-i\frac{1}{\sqrt{\gamma}}x}$ and
$\lambda=1$.
\item[(iii)] For any fixed $\epsilon>0$,
${\cal V}_{-\frac{1}{\sqrt{\gamma}}+\epsilon,\,0}(q/p') <\infty$. In
fact,
$$
\frac{q}{p'} = \frac{(1+\gamma \xi^2)e^{ix\xi}}{\alpha(1-\gamma \xi^2)},\quad
{\cal V}_{-\frac{1}{\sqrt{\gamma}}+\epsilon,\,0}\Big(\frac{q}{p'}\Big) =
\int_{-\frac{1}{\sqrt{\gamma}}+\epsilon}^{0}
\left|\Big(\frac{q}{p'}\Big)'\right| d \xi < \infty.
$$
\item[(iv)] As $\xi \to 0-$, $q/p'  \to 1/\alpha$.
\end{enumerate}
Theorem \ref{phase} then implies
\begin{eqnarray}
\int_{-\frac{1}{\sqrt{\gamma}}}^0\frac{1}{1+\gamma \xi^2}\,e^{ix\xi
-i\frac{\alpha\xi}{1+\gamma \xi^2} t}\, d\xi
& \sim & e^{-\frac{\pi i }{4}} \frac{1}{4}
e^{-i\frac{x}{\sqrt{\gamma}}}\,\Gamma\left(\frac{1}{2}\right) \frac{e^{it
(\frac{\alpha}{2\sqrt{\gamma}})}}{(\alpha/2\sqrt{\gamma}t)^{1/2}}
\nonumber \\
&=&
\frac{\sqrt{\pi}\,e^{-\frac{\pi i }{4}}\,
e^{\frac{i}{\sqrt{\gamma}}(-x+\frac{\alpha t}{2})}}{2\,\sqrt{2\alpha}
\sqrt[4]{\gamma}\,
\sqrt{t}}. \nonumber
\end{eqnarray}

\vspace{.1in} The estimate (\ref{i4}) also follows from Theorem
\ref{phase}. The conditions can be similarly checked for this
integral. In fact, for
$$
a=\frac{1}{\sqrt{\gamma}},\quad b=M, \quad p(\xi)= -\frac{\alpha \xi}{1+\gamma \xi^2}
\quad\mbox{and}\quad q(\xi) =\frac{e^{i\xi x}}{1+\gamma \xi^2}.
$$
we have
$$
p(\xi)-p\Big(\frac{1}{\sqrt{\gamma}}\Big)  \,\,\sim \,\,
\frac{\alpha\sqrt{\gamma}}{2}
\Big(\xi-\frac{1}{\sqrt{\gamma}}\Big)^2,\qquad
q\Big(\frac{1}{\sqrt{\gamma}}\Big)\,\,\sim \,\,\frac12
e^{i\frac{1}{\sqrt{\gamma}}x}.
$$
It is also easy to verify that
${\cal V}_{\frac{1}{\sqrt{\gamma}}+\epsilon,\,M}(q/p') <\infty$
for any fixed $\epsilon>0$. In
fact,
\begin{equation}\label{var}
\frac{q}{p'} = -\frac{(1+\gamma \xi^2)e^{ix\xi}}{\alpha(1-\gamma \xi^2)},\quad
{\cal V}_{\frac{1}{\sqrt{\gamma}}+\epsilon,\,M}\Big(\frac{q}{p'}\Big) =
\int_{\frac{1}{\sqrt{\gamma}}+\epsilon}^{M}
\left|\Big(\frac{q}{p'}\Big)'\right| d \xi < \infty.
\end{equation}
In addition, as $\xi \to M$, $q/p'$ tends to a finite limit. It then follows
from Theorem \ref{phase} that
$$
\int_{\frac{1}{\sqrt{\gamma}}}^M\frac{1}{1+\gamma \xi^2}\,e^{ix\xi
-i\frac{\alpha\xi}{1+\gamma \xi^2} t}\, d\xi
\quad \sim \quad
\frac{\sqrt{\pi}\,e^{\frac{\pi i }{4}}\,
e^{\frac{i}{\sqrt{\gamma}}(x-\frac{\alpha t}{2})}}{2\,\sqrt{2\alpha}
\sqrt[4]{\gamma}\,
\sqrt{t}}.
$$
When $M=\infty$, (\ref{var}) can not be verified and Theorem \ref{phase}
does not apply to the integral on $(-\infty,\infty)$.

\end{document}